\documentclass{amsart}   
\usepackage{amssymb}
\usepackage{amsmath}
\usepackage{amsfonts}

\theoremstyle{plain}

\numberwithin{equation}{section}

\begin{document}
\

%\vspace{7cm}
\begin{center}
\huge{\textbf{  Nonexistence results of global solutions \\
for fractional order
integral equations \\
on the Heisenberg group}} 
\end{center}

\vspace{0.5cm}

\begin{center}
\textbf{Abd Elhakim Lamairia\\
Department of mathematic and informatic. \\
university of Tebessa-Algeria \\
Email:hakim24039@gmail.com \\
abdelhakim.lamairia@univ-tebessa.dz}
\end{center}

\vspace{1cm}

\textbf{Abstract—}
\begin{normalsize}
We consider the fractional order integral equation with a time nonlocal nonlinearity
\[ ^{c}\mathbf{D}_{0\mid t}^{\beta}\left( u \right) +\left(-\Delta_{\mathbb{H}} 
\right)^{m} \left( u \right) = \frac{1}{\Gamma(\alpha)}\int_{0}^{t}\left( t-\omega\right) ^{\alpha-1}\vert u(\omega)\vert^{p} d\omega, \]
posed in $ (.,t)\in\mathbb{H}\times(0,\infty) $, supplemented with an initial data  $ u(.,0)=u_{0}(.) $,where $ m>1 \ , \ p>1 \ , \ 0<\beta<1 \ , \ 0<\alpha<1 $, and 
$ ^{c}\mathbf{D}_{0\mid t}^{\beta} $ denotes the  caputo fractional derivative of order $ \beta $, and 
$ \Delta_{\mathbb{H}} $ is the Laplacian operator on the $ (2N+1) $-dimensional Heisenberg group $ \mathbb{H} $.Then, we prove a blow up result for its solutions.  

\end{normalsize}

\vspace{0.25cm}

%\newpage 

\section{\large{\textbf{Introduction}}}

\vspace{0.7cm}

In this paper, we investigate the higher-order semilinear parabolic equation with
nonlocal in time nonlinearity of the following form:

\begin{equation*}
\left\lbrace \begin{array}{c}
\displaystyle ^{c}\mathbf{D}_{0\mid t}^{\beta}\left( u \right) +\left(-\Delta_{\mathbb{H}} 
\right)^{m} \left( u \right) = \mathbf{I}_{0\mid t}^{\alpha}\vert u(t)\vert^{p},   \\
  \eta=(x,y,\tau)\in \mathbb{H}, \ t>0 
\end{array} \right. \hspace{0.5cm} (1.1)
\end{equation*}\\

subject to the initial data
\[ u\left( \eta , 0\right)=u_{0}\left( \eta\right), \] 
Where $ \mathbf{I}_{0\mid t}^{\alpha}\psi $ is the Riemann–Liouville fractional integral of order  $(0 < \alpha < 1)$ defined for a continuous function  $\psi(t), t > 0$, 
\[(\mathbf{I}_{0\mid t}^{\alpha}\psi)(t)=\frac{1}{\Gamma(\alpha)}\int_{0}^{t}\left( t-\omega\right) ^{\alpha-1}\psi(\omega) d\omega, \]
Here, $ \Gamma(.) $ stands for the gamma function.\\
First, for the sake of the reader, we give some known facts about the Heisenberg
group $ \mathbb{H} $ and the operator $ \Delta_{\mathbb{H}} $. For their proof and more information, we refer
for example to  $[1, 4, 5, 11, 19]$ . The Heisenberg group $ \mathbb{H} $ , whose elements are $ \eta=(x,y,\tau) $ is the Lie group $ (\mathbb{R}^{2N+1},\circ) $ with the group operation $ "\circ" $ defined by
\[\eta\circ \tilde{\eta}=\left( x+\tilde{x},y+\tilde{y},\tau+\tilde{\tau}+2(<x,\tilde{y}>-<\tilde{x},y>)\right) ,\]
where $ <.,.> $ is the usual inner product in $ \mathbb{R}^{N} $, \ The laplacian $ \Delta_{\mathbb{H}} $ over $ \mathbb{H} $ is obtained from the vector fields $ X_{i}=\partial_{x_{i}}+2y_{i}\partial_{\tau} $ and $ Y_{i}=\partial_{y_{i}}+2x_{i}\partial_{\tau} $, by
\[ \Delta_{\mathbb{H}}=\sum_{i=1}^{N}\left( X_{i}^{2}+Y_{i}^{2}\right)  , \]
explicitly, we have
\[ \Delta_{\mathbb{H}}=\sum_{i=1}^{N}\left( \frac{\partial^{2}}{\partial x_{i}^{2}}+\frac{\partial^{2}}{\partial y_{i}^{2}}+4y_{i}\frac{\partial^{2}}{\partial x_{i}\partial \tau}-4x_{i}\frac{\partial^{2}}{\partial y_{i}\partial \tau}+4(x_{i}^{2}+y_{i}^{2})\frac{\partial^{2}}{\partial\tau^{2}}\right) , \]
A natural group of dilitations on $ \mathbb{H} $ is given by
\[ \delta_{\gamma}(\eta)=\left(\gamma x,\gamma y,\gamma^{2} \tau \right), \ \gamma>0,  \]
whose Jacobian determinant is $ \gamma^{Q} $
where \[Q=2N+2\]
is the homogeneous dimension of $ \mathbb{H} $.\\
The operator $ \Delta_{\mathbb{H}} $ is a degenerate elliptic operator. It is invariant with respect to the left translation of $ \mathbb{H} $ and homogeneous with respect to the dilatations $ \delta_{\gamma} $. More precisely, we have
\[ \Delta_{\mathbb{H}}\left( u(\eta \circ \tilde{\eta})\right) =\left(\Delta_{\mathbb{H}}u \right)(\eta \circ \tilde{\eta}), \ \Delta_{\mathbb{H}}(u\circ \delta_{\gamma})=\gamma^{2}(\Delta_{\mathbb{H}}u)\circ \delta_{\gamma} \ \ \eta , \tilde{\eta}\in \mathbb{H}.   \]
The natural distance from $ \eta $ to the origin is
\[ \vert\eta\vert_{\mathbb{H}}=\left(\tau^{2}+\left(\sum_{i=1}^{N}\left( x_{i}^{2}+y_{i}^{2}\right)  \right)^{2}  \right)^{\frac{1}{4}}.  \]
Now,we call sub-elliptic gradient
\[ \nabla_{\mathbb{H}}=\left( X,Y\right) =\left( X_{1},...,X_{N},Y_{1},...,Y_{N}\right) , \]
A remarkable property of the Kohn Laplacian is that a fundamental solution of $ -\Delta_{\mathbb{H}} $ with pole at zero is given by \[ \Gamma(\eta)=\frac{C_{\Lambda}}{ \vert \eta\vert_{\mathbb{H}}^{\Lambda-2} }, \]
where $ C_{\Lambda} $ is a suitable positive constant.\\
A basic role in the functional analysis on the Heisenberg group is played by
the following Sobolev-type inequality
\[ \Vert v\Vert_{\Lambda^{*}}^{2}=c\Vert \nabla_{\mathbb{H}}v\Vert_{2}^{2}, \forall v\in C_{0}^{\infty}(\mathbb{H}^{N}),  \]
where $ \Lambda^{*}=\frac{2\Lambda}{\Lambda-2} $ and $ c $ is a positive constant.\\
This inequality ensures in particular that for every domain $ \Omega $ the function
\[ \Vert v\Vert\leq \Vert \nabla_{\mathbb{H}}v\Vert_{2}, \]

is a norm on $ C_{0}^{\infty}(\Omega) $. We denote by $ S_{0}^{1}(\Omega) $ the closure of $ C_{0}^{\infty}(\Omega) $ with respect to
this norm; $ S_{0}^{1}(\Omega) $ becomes a Hilbert space with the inner product
\[ <u,v>_{S_{0}^{1}}=\int_{\Omega}<\nabla_{\mathbb{H}}u,\nabla_{\mathbb{H}}v>,
 \]
\textbf{Fractional powers of sub-elliptic Laplacians.} Here, we recall a result on
fractional powers of sub-Laplacian in the Heisenberg group. Let $ \mathit{N}(t,x) $ be the fundamental solution of $ \Delta_{\mathbb{H}}+\frac{\partial}{\partial t} $ . For all $ 0<\beta<4 $, the integral
\[ \mathit{R}_{\beta}(x)=\frac{1}{\Gamma\left(\frac{\beta}{2} \right) }\int_{0}^{+\infty}t^{\frac{\beta}{2}-1}\mathit{N}(t,x)dt, \]
converges absolutely for $ x\neq 0 $. If $ \beta<0 , \beta \neq 0,-2,-4,... $,then

\[ \tilde{\mathit{R}}_{\beta}(x)=\frac{\frac{\beta}{2}}{\Gamma\left(\frac{\beta}{2} \right) }\int_{0}^{+\infty}t^{\frac{\beta}{2}-1}\mathit{N}(t,x)dt, \]
defines a smooth function in $ \mathbb{H} -\left\lbrace 0\right\rbrace  $, since $ t\mapsto \mathit{N}(t,x) $, vanishes of infinite order as $ t\rightarrow 0 $ if $ x\neq 0 $. In addition, $ \tilde{\mathit{R}}_{\beta} $ is positive and $ \mathbb{H} $-homogeneous of degree $ \beta-4 $.\\
\\
\textbf{Theorem:}\\
For every $ v\in S(\mathbb{H}) $ (Schwartz’s class), we have $ \left(-\Delta_{\mathbb{H}} \right)^{s} \in L^{2}(\mathbb{H})  $ and
\[ \left(-\Delta_{\mathbb{H}} \right)^{s}=\int_{\mathbb{H}}\left( v(x \circ y)-v(x)-\chi(y)<\nabla_{\mathbb{H}}v(x),y>\right)\tilde{\mathit{R}}_{-2s}(y)dy,  \] 

where $ \chi $ is the characteristic function of the unit ball $ \mathit{B}_{\rho}(0,1) $, ($ \rho(x)=\mathit{R}_{2-\alpha}^{\frac{-1}{2-\alpha}}(x),$  \ \ \ \ \ \  $0<\alpha<2 $ , $ \rho $ is an $ \mathbb{H} $-homogeneous norm in $ \mathbb{H} $ smooth outside the origin).\\

\vspace{0.5cm}
\section{\textbf{Preliminaries}}

\vspace{0.5cm}

\subsection{\textbf{Definition}}

(Riemann-Liouville fractional derivatives)

Let $ f\in AC[a , b] , -\infty < a<b <+\infty  $,\footnote{let $ AC[a, b] $ be the space of functions $ f $ which are absolutely continuous on $[a,b]$.\\
$ AC^{n}\left[ a,b\right] =\left\lbrace f:\left[ a,b\right]\rightarrow \mathbb{C} \  \mbox{and}  \left(D^{n-1}f \right)(x)\in AC[a, b] \ \ \left( D=\frac{d}{dx}\right) \right\rbrace   $.\\
In particular, $ AC^{1}\left[ a,b\right]=AC[a, b], $ } The Riemann-Liouville left- and right-sided
fractional derivatives of order $ \alpha \in (0 , 1) $ are, respectively, defined by

\[ \displaystyle \mathbf{D}_{a\vert t}^{\alpha}f(t) = \frac{d}{dt}\mathbf{I}_{a\vert t}^{1-\alpha}f(t)\]

\[=\frac{1}{\Gamma \left( 1-\alpha \right) }\frac{d}{dt}\int_{a}^{t}
(t-\tau)^{-\alpha}f(\tau)d\tau, \ \ t>a  \hspace{1cm} (2.1) \]

and

\[ \displaystyle \mathbf{D}_{t\vert b}^{\alpha}f(t) = -\frac{d}{dt}\mathbf{I}_{t\vert b}^{1-\alpha}f(t) \]

\[=-\frac{1}{\Gamma \left( 1-\alpha \right) }\frac{d}{dt}\int_{t}^{b}
(\tau-t)^{-\alpha}f(\tau)d\tau, \ \ t<b  \hspace{1cm} (2.2) \]

\subsection{\textbf{Definition}}

(Riemann-Liouville fractional integrals)

Let $ f\in L^{1}(a , b) ,-\infty < a<b <+\infty  $, The Riemann-Liouville left- and right-sided
fractional integrals of order $ \alpha \in (0 , 1) $ are, respectively, defined by 

\[ \mathbf{I}_{a\mid t}^{\alpha}f(t)=\frac{1}{\Gamma(\alpha)}\int_{a}^{t}(t-\tau)^{-(1-\alpha)}f(\tau)d\tau , \ \ t>a  \hspace{1cm}  (2.3) \]
and
\[ \mathbf{I}_{t\mid b}^{\alpha}f(t)=\frac{1}{\Gamma(\alpha)}\int_{t}^{b}(\tau-t)^{-(1-\alpha)}f(\tau)d\tau , \ \ t<b  \hspace{1cm}  (2.4) \]

\subsection{\textbf{Definition}}

For $ 0<\alpha<1 $, the Caputo derivative of order $ \alpha $ for a differentiable function $ f : [0, \infty) \rightarrow  \mathbb{R} $ can been written as
\[ \displaystyle ^{c}\mathbf{D}_{a\vert t}^{\alpha}f(t) =\frac{1}{\Gamma \left( 1-\alpha \right) }\frac{d}{dt}\int_{a}^{t}
(t-\tau)^{-\alpha}f^{'}(\tau)d\tau, \ \ t>a  \hspace{1cm} (2.5) \]
It is clear that
\[ ^{c}\mathbf{D}_{a\vert t}^{\alpha}f(t) =\mathbf{D}_{a\vert t}^{\alpha}\left[ f(t)-f(0)\right] , \]

Finally, taking into account the following integration by parts
formula:
\begin{equation*}
\int_{a}^{b}f\left( t\right) \mathbf{D}_{a\mid t}^{\alpha }g \left(
t\right)dt =\int_{a}^{b} \mathbf{D}_{t\mid b}^{\alpha }f \left( t\right)
g\left( t\right)dt.
\end{equation*}

\subsection{\textbf{Proposition}}

For $ 0<\alpha<1 ,-\infty < a<b <+\infty  $, we have the following identities
\[ \displaystyle \mathbf{D}_{a\vert t}^{\alpha}\mathbf{I}_{a\mid t}^{\alpha}f(t)=f(t) , \ t\in (a , b) \]	
for all $ f\in L^{r}(a , b) , 1\leq r \leq \infty $\\
and
\[ -\mathbf{D}\mathbf{D}_{t\mid b}^{\alpha}f=\mathbf{D}_{t\mid b}^{1+\alpha}f , \]
for all $ f\in AC^{2}[a , b] , $ where $ \mathbf{D}=\frac{d}{dt} $.\\

For $ \rho\gg 1 $ and $ 0<\alpha<1 $. \ Let

\begin{equation*}
f(t)=\left\lbrace \begin{array}{c}
\displaystyle \left( 1-\frac{t}{T}\right)^{\rho} , \ \ \ 0< t \leq T, \\
  0, \ \ \ \ \ \ \ t\geq T, 
\end{array} \right. \hspace{0.5cm} (2.6)
\end{equation*}

\[ \mathbf{D}_{t\mid T}^{\alpha}f(t)=\frac{(1-\alpha+\rho)\Gamma( \rho-1)}{\Gamma(2-\alpha-\rho)}T^{-\alpha}\left( 1-\frac{t}{T}\right)^{\rho-\alpha}, \]

and
\[ \int_{0}^{T}f(t)^{\frac{-p^{'}}{p}}\vert\mathbf{D}_{t\mid T}^{\alpha}f(t)\vert^{p^{'}}=CT^{1-p^{'}\alpha}, \]

\section{\textbf{ Nonexistence results}}

\vspace{0.3cm}

\subsection{\textbf{Definition}}

(Weak solution).Let $ T>0 $, a locally integrable function $ u\in C\left( [0,T], L_{loc}^{1}\left( Q_{T}\right) \cap L_{loc}^{p}\left( Q_{T}\right)\right)  $ is called a local weak
solution of (1.1) in $ Q_{T} $ $ \left( Q_{T}=\mathbb{H}\times [0,T]\right) $
subject to the initial data
$ u_{0}\in L_{loc}^{1}\left( \mathbb{H}\right)  $ if the equality
\[ \int_{Q_{T}}u_{0} \mathbf{D}_{t\mid T}^{\beta} \varphi d\omega + \int_{Q_{T}}\varphi \mathbf{I}_{0 \mid t}^{\alpha} \vert u \vert^{p}d\omega\\
 =\int_{Q_{T}}u\mathbf{D}_{t\mid T}^{\beta}\varphi d\omega +\int_{Q_{T}}u(-\Delta_{\mathbb{H}})^{m}\varphi d\omega \],
is satisfied for any $ \varphi $ be a smooth test function $ \varphi\in C_{0}^{\infty}(Q_{T})) $ witn
\[ \varphi(.,T)=0, \ \ \varphi\geq 0, \ \ d\omega=d\eta dt \]
and the solution is called global if $ T=+\infty. $\\

\subsection{\textbf{Theorem}}

Let $ p>1 $, and \\
\[  p < p_{c}=\frac{(2N+2)\beta+2m}{(2N+2)\beta+2m(1-\alpha)}, \] \ \ (c for critical)\\
Then, (1.1) does not have a nontrivial global weak solution.

%\vspace{5cm}
\subsection{\textbf{Proposition}}
Consider a convex function $ F\in C^{2}(\mathbb{R})$. Assume that $ \varphi\in C_{0}^{\infty}(\mathbb{R}^{2N+1}), $ then
\[ F^{'}(\varphi)(-\Delta_{\mathbb{H}})^{m}\varphi\geq (-\Delta_{\mathbb{H}})^{m}F(\varphi), \]
In particular, if $ F(0)=0 $ and $ \varphi\in C_{0}^{\infty}(\mathbb{R}^{2N+1}), $ then 
\[\int_{\mathbb{R}^{2N+1}}F^{'}(\varphi)(-\Delta_{\mathbb{H}})^{m}\varphi d\eta \geq 0. \]
Let us mention that hereafter we will use inequality (2.1) for $ F(\varphi)=\varphi^{l} , \ \ l\gg 1, \ \ \varphi\geq 0,$ in this case it reads

\begin{equation}
l \varphi^{l-1}(-\Delta_{\mathbb{H}})^{m}\varphi\geq (-\Delta_{\mathbb{H}})^{m}\varphi^{l},
\end{equation}
We need the following Lemma taken from [32].

\subsection{\textbf{Lemma}}
Let $ f\in L^{1}(\mathbb{R}^{2N+1}) $ and $ \int_{\mathbb{R}^{2N+1}}f d\eta\geq 0. $
Then there exists a test function $ 0\leq \varphi \leq 1, $ such that
\begin{equation}
 \int_{\mathbb{R}^{2N+1}}f\varphi d\eta\geq 0.
\end{equation}
\\  
\textbf{Proof of theorem :}.\\

The proof is done by contradiction. Suppose that $ u $ is a global bounded weak solution. First we Choose the test function.  For this aim, we shall use a non-negative smooth function $ \phi $ which was constructed in [20].

\begin{equation}
\phi(\xi)=\left\lbrace \begin{array}{c}
 1 \ \ \ \ if \ \ 0\leq \xi \leq 1.\\
\searrow \ \ \ \ if \ \ 1\leq \xi \leq 2,\\
0 \ \ \ \ if \ \ \ \ \ \xi \geq 2,
\end{array}\right. 
\end{equation}

\[ \varphi_{1}(\eta)=\phi\left( \frac{\tau^{2}+\vert x\vert^{4}+\vert y\vert^{4}}{R^{4}}\right) , \  \ \ \eta=(x,y,\tau)\in \mathbb{H} \]

\[ \varphi_{2}(t)=\left\lbrace \begin{array}{c}
\displaystyle \left( 1-\frac{t}{T}\right)^{\rho} , \ \ \ 0< t \leq T, \\
  0, \ \ \ \ \ \ \ t\geq T, 
\end{array} \right. \ \ \rho\gg 1 \]

\[ \varphi(\eta,t)=\mathbf{D}_{t\mid TR^{\frac{2m}{\beta}}}^{\alpha}\tilde{\varphi}(\eta,t)=\varphi_{1}^{l}\left( \eta\right) \mathbf{D}_{t\mid TR^{\frac{2m}{\beta}}}^{\alpha}\varphi_{2}\left(\frac{t}{R^{\frac{2m}{\beta}}} \right) ,  \ \ R>0  \]
Let, $ Q=\mathbb{H}\times \left[ 0,TR^{\frac{2m}{\beta}}\right],  $ \\

Using the Definition 3.1, we obtain

\[ \int_{Q}u_{0} \mathbf{D}_{t\mid TR^{\frac{2m}{\beta}}}^{\beta} \mathbf{D}_{t\mid TR^{\frac{2m}{\beta}}}^{\alpha}\tilde{\varphi}(\eta,t)d\eta dt + \int_{Q}\mathbf{D}_{t\mid TR^{\frac{2m}{\beta}}}^{\alpha}\tilde{\varphi}(\eta,t) \mathbf{I}_{0 \mid t}^{\alpha} \vert u \vert^{p}d\eta dt \]
\[ =\int_{Q}u\mathbf{D}_{t\mid TR^{\frac{2m}{\beta}}}^{\beta}\mathbf{D}_{t\mid TR^{\frac{2m}{\beta}}}^{\alpha}\tilde{\varphi}(\eta,t) d\eta dt +\int_{Q}u(-\Delta_{\mathbb{H}})^{m}\mathbf{D}_{t\mid TR^{\frac{2m}{\beta}}}^{\alpha}\tilde{\varphi}(\eta,t) d\eta dt, \]
A simple computation yields $ \mathbf{D}_{t\mid TR^{\frac{2m}{\beta}}}^{\beta} \left( \mathbf{D}_{t\mid TR^{\frac{2m}{\beta}}}^{\alpha}\tilde{\varphi}\right)= \mathbf{D}_{t\mid TR^{\frac{2m}{\beta}}}^{\alpha+\beta}\tilde{\varphi},  $
we obtain
\[ c(TR^{\frac{2m}{\beta}})^{1-(\alpha+\beta)}\int_{\mathbb{H}}u_{0} \varphi_{1}^{l}(\eta)d\eta + \int_{Q}\tilde{\varphi}\vert u \vert^{p}d\eta dt \]

\[ =\int_{Q}u\varphi_{1}^{l}(\eta)\mathbf{D}_{t\mid TR^{\frac{2m}{\beta}}}^{\alpha+\beta}\varphi_{2}(\frac{t}{R^{\frac{2m}{\beta}}})d\eta dt + \int_{Q}u(-\Delta_{\mathbb{H}})^{m}\varphi_{1}^{l}(\eta)\mathbf{D}_{t\mid TR^{\frac{2m}{\beta}}}^{\alpha}\varphi_{2}(\frac{t}{R^{\frac{2m}{\beta}}})d\eta dt,   \]

The application of inequality (3.1)
\[ l \varphi_{1}^{l-1}(-\Delta_{\mathbb{H}})^{m}\varphi_{1} \geq (-\Delta_{\mathbb{H}})^{m}\varphi_{1}^{l}, \]
implies that 
\[ \int_{Q}\vert u \vert^{p}\tilde{\varphi} d\eta dt \]
\[ \leq l\int_{Q}u\varphi_{1}^{l-1}(\eta)(-\Delta_{\mathbb{H}})^{m}\varphi_{1}(\eta)\mathbf{D}_{t\mid TR^{\frac{2m}{\beta}}}^{\alpha}\varphi_{2}(\frac{t}{R^{\frac{2m}{\beta}}})d\eta dt + \int_{Q}u\varphi_{1}^{l}(\eta)\mathbf{D}_{t\mid TR^{\frac{2m}{\beta}}}^{\alpha+\beta}\varphi_{2}(\frac{t}{R^{\frac{2m}{\beta}}})d\eta dt,\]

For estimating the second member of the above inequality, we write
\[ \int_{Q}u\varphi_{1}^{l-1}(\eta)(-\Delta_{\mathbb{H}})^{m}\varphi_{1}(\eta)\mathbf{D}_{t\mid TR^{\frac{2m}{\beta}}}^{\alpha}\varphi_{2}(\frac{t}{R^{\frac{2m}{\beta}}})d\eta dt \]

\[=\int_{Q}u \tilde{\varphi}^{\frac{1}{p}} \varphi_{1}^{l-1}(\eta)(-\Delta_{\mathbb{H}})^{m}\varphi_{1}(\eta)\mathbf{D}_{t\mid TR^{\frac{2m}{\beta}}}^{\alpha}\varphi_{2}(\frac{t}{R^{\frac{2m}{\beta}}})  \tilde{\varphi}^{\frac{-1}{p}}d\eta dt.  \]

According to $ \epsilon $-Young inequality
\[XY \leq \epsilon X^{p}+C(\epsilon)Y^{p^{'}}, \ \ \ p+p^{'}=pp^{'},    \]
we have

\[ \int_{Q}u\varphi_{1}^{l-1}(\eta)(-\Delta_{\mathbb{H}})^{m}\varphi_{1}(\eta)\mathbf{D}_{t\mid TR^{\frac{2m}{\beta}}}^{\alpha}\varphi_{2}(\frac{t}{R^{\frac{2m}{\beta}}})d\eta dt  \]

\[ \leq \epsilon \int_{Q}\vert u\vert^{p}\tilde{\varphi}d\eta dt + 
C_{1}(\epsilon)\int_{Q}\varphi_{1}^{(l-1)p^{'}}(\eta)\vert(-\Delta_{\mathbb{H}})^{m}\varphi_{1}(\eta)\mathbf{D}_{t\mid TR^{\frac{2m}{\beta}}}^{\alpha}\varphi_{2}(\frac{t}{R^{\frac{2m}{\beta}}})\vert^{p^{'}}  \tilde{\varphi}^{\frac{-p^{'}}{p}}d\eta dt. \]
In the same way, we get

\[ \int_{Q}u\varphi_{1}^{l}(\eta)\mathbf{D}_{t\mid TR^{\frac{2m}{\beta}}}^{\alpha+\beta}\varphi_{2}(\frac{t}{R^{\frac{2m}{\beta}}})d\eta dt \]

\[ \leq \epsilon \int_{Q}\vert u\vert^{p}\tilde{\varphi}d\eta dt +
C_{2}(\epsilon)\int_{Q}\vert\varphi_{1}^{l}(\eta)\mathbf{D}_{t\mid TR^{\frac{2m}{\beta}}}^{\alpha+\beta}\varphi_{2}(\frac{t}{R^{\frac{2m}{\beta}}})\vert^{p^{'}}\tilde{\varphi}^{\frac{-p^{'}}{p}} d\eta dt  \]

Now, when $ \epsilon $ is small, and $ C=\max \left\lbrace C_{1}(\epsilon) , C_{2}(\epsilon)\right\rbrace   $ we obtain

\[ \int_{Q}\vert u \vert^{p}\tilde{\varphi} d\eta dt \leq C\left\lbrace \int_{Q}\varphi_{1}^{(l-1)p^{'}}(\eta)\vert(-\Delta_{\mathbb{H}})^{m}\varphi_{1}(\eta)\mathbf{D}_{t\mid TR^{\frac{2m}{\beta}}}^{\alpha}\varphi_{2}(\frac{t}{R^{\frac{2m}{\beta}}})\vert^{p^{'}}  \tilde{\varphi}^{\frac{-p^{'}}{p}}d\eta dt\right.   \]
 
\[\left.  + \int_{Q}\vert\varphi_{1}^{l}(\eta)\mathbf{D}_{t\mid TR^{\frac{2m}{\beta}}}^{\alpha+\beta}\varphi_{2}(\frac{t}{R^{\frac{2m}{\beta}}})\vert^{p^{'}}\tilde{\varphi}^{\frac{-p^{'}}{p}} d\eta dt \right\rbrace,   \]
as
\[ \tilde{\varphi}^{\frac{-p^{'}}{p}}(\eta,t)=\varphi_{1}^{\frac{-p^{'}}{p}l}(\eta)\varphi_{2}^{\frac{-p^{'}}{p}}(\frac{t}{R^{\frac{2m}{\beta}}}), \ \ \ p^{'}=\frac{p}{p-1}  \]

we have
\[ \int_{Q}\vert u \vert^{p}\tilde{\varphi} d\eta dt \leq C\left\lbrace \int_{Q}\varphi_{1}^{(l-p^{'})}(\eta)\varphi_{2}^{\frac{-1}{p-1}}(\frac{t}{R^{\frac{2m}{\beta}}})\vert(-\Delta_{\mathbb{H}})^{m}\varphi_{1}(\eta)\mathbf{D}_{t\mid TR^{\frac{2m}{\beta}}}^{\alpha}\varphi_{2}(\frac{t}{R^{\frac{2m}{\beta}}})\vert^{p^{'}} d\eta dt \right. \] 

\[ \left. + \int_{Q}\varphi_{1}^{l}(\eta)\varphi_{2}^{\frac{-1}{p-1}}(\frac{t}{R^{\frac{2m}{\beta}}})\vert\mathbf{D}_{t\mid TR^{\frac{2m}{\beta}}}^{\alpha+\beta}\varphi_{2}(\frac{t}{R^{\frac{2m}{\beta}}})\vert^{p^{'}} d\eta dt \right\rbrace,   \]

We apply the change of next variables $ \tilde{\tau}=\frac{\tau}{R^{2}} $, \ \ $ \tilde{x}=\frac{x}{R} $, \ \ $ \tilde{y}=\frac{y}{R} $, \ \ $ \tilde{t}=\frac{t}{R^{\frac{2m}{\beta}}} $, then we put  \[\Omega=\left\lbrace \tilde{\eta}=(\tilde{x},\tilde{y},\tilde{\tau})\in \mathbb{H};\ 0\leq \tilde{\tau}^{2}+\vert \tilde{x}\vert^{4}+ \vert \tilde{y}\vert^{4} \leq 2\right\rbrace \].

if we put
\[ \mathcal{A}= \int_{Q}\varphi_{1}^{(l-p^{'})}(\eta)\varphi_{2}^{\frac{-1}{p-1}}(\frac{t}{R^{\frac{2m}{\beta}}})\vert(-\Delta_{\mathbb{H}})^{m}\varphi_{1}(\eta)\mathbf{D}_{t\mid TR^{\frac{2m}{\beta}}}^{\alpha}\varphi_{2}(\frac{t}{R^{\frac{2m}{\beta}}})\vert^{p^{'}} d\eta dt, \]

\[ \mathcal{B}=\int_{Q}\varphi_{1}^{l}(\eta)\varphi_{2}^{\frac{-1}{p-1}}(\frac{t}{R^{\frac{2m}{\beta}}})\vert\mathbf{D}_{t\mid TR^{\frac{2m}{\beta}}}^{\alpha+\beta}\varphi_{2}(\frac{t}{R^{\frac{2m}{\beta}}})\vert^{p^{'}} d\eta dt, \]

we get
\[ \mathcal{A}=\int_{0}^{TR^{\frac{2m}{\beta}}}\varphi_{2}^{\frac{-1}{p-1}}(\frac{t}{R^{\frac{2m}{\beta}}})\vert \mathbf{D}_{t\mid TR^{\frac{2m}{\beta}}}^{\alpha}\varphi_{2}(\frac{t}{R^{\frac{2m}{\beta}}})\vert^{p^{'}}dt \int_{\mathbb{H}} \varphi_{1}^{(l-p^{'})}(\eta)\vert(-\Delta_{\mathbb{H}})^{m}\varphi_{1}(\eta)\vert^{p^{'}} d\eta  \]

\[ =R^{\frac{2m}{\beta}-\frac{2mp\alpha}{\beta(p-1)}} \int_{0}^{T}\varphi_{2}^{\frac{-p^{'}}{p}}\left( \tilde{t}\right) \vert \mathbf{D}_{\tilde{t}\mid T}^{\alpha}\varphi_{2}\left( \tilde{t}\right) \vert^{p^{'}} d\tilde{t} \] 

\[\times R^{2N+2}\int_{\Omega} \phi^{(l-p^{'})}\left( \tilde{\tau}^{2}+\vert \tilde{x}\vert^{4}+ \vert \tilde{y}\vert^{4}\right) \vert(-\Delta_{\mathbb{H}})^{m}\phi\left( \tilde{\tau}^{2}+\vert \tilde{x}\vert^{4}+ \vert \tilde{y}\vert^{4}\right) \vert^{p^{'}} d\tilde{x} d\tilde{y} d\tilde{\tau},  \]

\[ =CT^{1-\frac{p\alpha}{p-1}}R^{2N+2+\frac{2m}{\beta}-\frac{2mp\alpha}{\beta(p-1)}}\] 

\[\times \int_{\Omega} \phi^{(l-p^{'})}\left( \tilde{\tau}^{2}+\vert \tilde{x}\vert^{4}+ \vert \tilde{y}\vert^{4}\right) \vert(-\Delta_{\mathbb{H}})^{m}\phi\left( \tilde{\tau}^{2}+\vert \tilde{x}\vert^{4}+ \vert \tilde{y}\vert^{4}\right) \vert^{p^{'}} d\tilde{x} d\tilde{y} d\tilde{\tau},  \]
and

\[ \mathcal{B}=\int_{0}^{TR^{\frac{2m}{\beta}}}\varphi_{2}^{\frac{-1}{p-1}}(\frac{t}{R^{\frac{2m}{\beta}}})\vert \mathbf{D}_{t\mid TR^{\frac{2m}{\beta}}}^{\alpha+\beta}\varphi_{2}(\frac{t}{R^{\frac{2m}{\beta}}})\vert^{p^{'}}dt \int_{\mathbb{H}} \varphi_{1}^{l}(\eta) d\eta  \]
\[ =R^{\frac{2m}{\beta}-\frac{2mp(\alpha+\beta)}{\beta(p-1)}} \int_{0}^{T}\varphi_{2}^{\frac{-p^{'}}{p}}\left( \tilde{t}\right) \vert \mathbf{D}_{\tilde{t}\mid T}^{\alpha+\beta}\varphi_{2}\left( \tilde{t}\right) \vert^{p^{'}} d\tilde{t} \] 

\[\times R^{2N+2}\int_{\Omega} \phi^{l}\left( \tilde{\tau}^{2}+\vert \tilde{x}\vert^{4}+ \vert \tilde{y}\vert^{4}\right) d\tilde{x} d\tilde{y} d\tilde{\tau},  \]

\[ =CT^{1-\frac{p(\alpha+\beta)}{p-1}}R^{2N+2+\frac{2m}{\beta}-\frac{2mp(\alpha+\beta)}{\beta(p-1)}} \times \int_{\Omega} \phi^{l}\left( \tilde{\tau}^{2}+\vert \tilde{x}\vert^{4}+ \vert \tilde{y}\vert^{4}\right) d\tilde{x} d\tilde{y} d\tilde{\tau},  \]
in the last 
\begin{equation}
 \int_{Q}\vert u \vert^{p}\tilde{\varphi} d\eta dt \leq C\left\lbrace \mathcal{A}+\mathcal{B}\right\rbrace \leq \mathcal{C}R^{2N+2+\frac{2m}{\beta}-\frac{2mp\alpha}{\beta(p-1)}} 
\end{equation}
Now, if
\[ 2N+2+\frac{2m}{\beta}-\frac{2mp\alpha}{\beta(p-1)}<0 \Leftrightarrow p<p_{c}  \]    

by letting $ R\rightarrow +\infty $ in (3.4), we obtain

\[ \int_{Q}\vert u \vert^{p}d\eta dt=0 \Rightarrow u\equiv 0, \]\\

this is a contradiction. \ \ \ \  $\square$

\end{document}